\documentclass[12pt,leqno]{smfart}

\usepackage[frenchb]{babel}
\usepackage{amsfonts}

\newdimen\DiagramCellHeight\DiagramCellHeight3em 
\newdimen\DiagramCellWidth\DiagramCellWidth3em 
\newdimen\MapBreadth\MapBreadth.04em 
\newdimen\MapShortFall\MapShortFall.4em 
\newdimen\PileSpacing\PileSpacing1em 
\def\labelstyle{\ifincommdiag\textstyle\else\scriptstyle\fi}
\let\objectstyle\displaystyle

\def\rArr{\HorizontalMap\empty-\empty-\rhla}

\font\tenln=line10

\mathchardef\lt="313C \mathchardef\gt="313E

\def\dhvee{\vbox\tozpt{\vss\hbox{$\vee$}\kern0pt}}
\def\dhcvee{\vbox\tozpt{\vss\hbox{$\curlyvee$}\kern0pt}}
\def\dhvvee{\vbox\tozpt{\vss\hbox{$\vee$}\kern-.6ex\hbox{$\vee$}\kern0pt}}
\def\uhvvee{\vbox\tozpt{\hbox{$\wedge$}\kern-.6ex\hbox{$\wedge$}\vss}}
\def\twoheaddownarrow{\rlap{$\downarrow$}\raise-.5ex\hbox{$\downarrow$}}
\def\twoheaduparrow{\rlap{$\uparrow$}\raise.5ex\hbox{$\uparrow$}}
\def\rhla{\vbox\tozpt{\vss\hbox\tozpt{\hss\tenln\char'55}\kern\axisheight}}
\def\lhla{\vbox\tozpt{\vss\hbox\tozpt{\tenln\char'33\hss}\kern\axisheight}}
\def\rthooka{\raise.603ex\hbox{$\scriptscriptstyle\subset$}}
\def\lthooka{\raise.603ex\hbox{$\scriptscriptstyle\supset$}}
\def\rthookb{\raise-.022ex\hbox{$\scriptscriptstyle\subset$}}
\def\lthookb{\raise-.022ex\hbox{$\scriptscriptstyle\supset$}}

\def\SEpbk{\rlap{\smash{\kern0.1em \vrule depth 2.67ex height -2.55ex width 0%
.9em \vrule height -0.46ex depth 2.67ex width .05em }}}
\def\SWpbk{\llap{\smash{\vrule height -0.46ex depth 2.67ex width .05em \vrule
depth 2.67ex height -2.55ex width .9em \kern0.1em }}}
\def\NEpbk{\rlap{\smash{\kern0.1em \vrule depth -3.48ex height 3.67ex width 0%
.95em \vrule height 3.67ex depth -1.39ex width .05em }}}
\def\NWpbk{\llap{\smash{\vrule height 3.6ex depth -1.39ex width .05em \vrule
depth -3.48ex height 3.67ex width .95em \kern0.1em }}}

\newcount\cdna\newcount\cdnb\newcount\cdnc\newcount\cdnd\cdna=\catcode`\@%
\catcode`\@=11 \let\then\relax\def\loopa#1\repeat{\def\bodya{#1}\iteratea}%
\def\iteratea{\bodya\let\next\iteratea\else\let\next\relax\fi\next}\def\loopb
#1\repeat{\def\bodyb{#1}\iterateb}\def\iterateb{\bodyb\let\next\iterateb\else
\let\next\relax\fi\next} \ifx\inputlineno\undefined\then\def\theinputlineno#1%
{}\let\noerrcontext\relax\else\def\noerrcontext{\errorcontextlines=-1}\def
\theinputlineno#1{ #1 line \number\inputlineno}\fi\def\commdiagname{%
Commutative Diagrams}\def\mapclasherr{\message{! \commdiagname: clashing maps%
\theinputlineno{in diagram ending at}}}\def\mapctxterr#1{\noerrcontext
\errhelp={This combination of maps doesn't make geometrical sense.}%
\expandafter\errmessage{\commdiagname: #1 \whyforbidmap is not allowed}}\def
\Obsmess#1#2{\expandafter\message{! \commdiagname\space Warning: #1 \string#2
is obsolete}\expandafter\message{ (ignored \theinputlineno{- first occurred on%
})}}\def\ObsDim#1{\Obsmess{Dimension}{#1}\global\let#1\ObsDimq\ObsDimq}\def
\ObsDimq{\dimen@=}\def\HorizontalMapLength{\ObsDim\HorizontalMapLength}\def
\VerticalMapHeight{\ObsDim\VerticalMapHeight}\def\VerticalMapDepth{\ObsDim
\VerticalMapDepth}\def\VerticalMapExtraHeight{\ObsDim\VerticalMapExtraHeight}%
\def\VerticalMapExtraDepth{\ObsDim\VerticalMapExtraDepth}\def\ObsCount#1{%
\Obsmess{Count}{#1}\global\let#1\ObsCountq\ObsCountq}\def\ObsCountq{\count@=}%
\def\DiagonalLineSegments{\ObsCount\DiagonalLineSegments}\def\tozpt{to\z@}%
\def\sethorizhtdp{\dimen8=\axisheight\dimen9=\MapBreadth\advance\dimen8.5%
\dimen9\advance\dimen9-\dimen8}\def\horizhtdp{height\dimen8 depth\dimen9 }%
\def\axisheight{\fontdimen22\the\textfont2 }\countdef\boxc@unt=14

\def\bombparameters{\hsize\z@\rightskip\z@ plus1fil minus\maxdimen
\parfillskip\z@\linepenalty9000 \looseness0 \hfuzz\maxdimen\hbadness10000
\clubpenalty0 \widowpenalty0 \displaywidowpenalty0 \interlinepenalty0
\predisplaypenalty0 \postdisplaypenalty0 \interdisplaylinepenalty0
\interfootnotelinepenalty0 \floatingpenalty0 \brokenpenalty0 \everypar{}%
\leftskip\z@\parskip\z@\parindent\z@\pretolerance10000 \tolerance10000
\hyphenpenalty10000 \exhyphenpenalty10000 \binoppenalty10000 \relpenalty10000
\adjdemerits0 \doublehyphendemerits0 \finalhyphendemerits0 \prevdepth\z@}\def
\startbombverticallist{\hbox{}\penalty1\nointerlineskip}

\def\pushh#1\to#2{\setbox#2=\hbox{\box#1\unhbox#2}}\def\pusht#1\to#2{\setbox#%
2=\hbox{\unhbox#2\box#1}}

\newif\ifallowhorizmap\allowhorizmaptrue\newif\ifallowvertmap
\allowvertmapfalse\def\whyforbidmap{outside diagram }\newif\ifincommdiag
\incommdiagfalse

\def\diagram{\leavevmode\hbox\bgroup$\vcenter\bgroup\startbombverticallist
\incommdiagtrue\baselineskip\DiagramCellHeight\lineskip\z@\lineskiplimit\z@
\mathsurround\z@\tabskip\z@\let\\\diagcr\allowhorizmaptrue\allowvertmaptrue
\def\whyforbidmap{at the beginning of a cell }\halign\bgroup\lcdtempl##%
\rcdtempl&&\lcdtempl##\rcdtempl\cr}\def\enddiagram{\crcr\egroup
\reformatmatrix\egroup$\egroup}

\def\lcdtempl{\futurelet\thefirsttoken\dolcdtempl}\newif\ifemptycell\def
\dolcdtempl{\ifx\thefirsttoken\rcdtempl\then\hskip1sp plus 1fil \emptycelltrue
\else\hfil$\emptycellfalse\objectstyle\fi}\def\rcdtempl{\ifemptycell\else$%
\hfil\fi}\def\diagcr{\cr} \def\across#1{\span\omit\mscount=#1 \loop\ifnum
\mscount>2 \spAn\repeat\ignorespaces}\def\spAn{\relax\span\omit\advance
\mscount by -1}

\def\CellSize{\afterassignment\cdhttowd\DiagramCellHeight}\def\cdhttowd{%
\DiagramCellWidth\DiagramCellHeight}\def\MapsAbut{\MapShortFall\z@}

\newcount\cdvdl\newcount\cdvdr\newcount\cdvd\newcount\cdbfb\newcount\cdbfr
\newcount\cdbfl\newcount\cdvdr\newcount\cdvdl\newcount\cdvd

\def\reformatmatrix{\bombparameters\cdvdl=\insc@unt\cdvdr=\cdvdl\cdbfb=%
\boxc@unt\advance\cdbfb1 \cdbfr=\cdbfb\setbox1=\vbox{}\dimen2=\z@\loop\setbox
0=\lastbox\ifhbox0 \dimen1=\lastskip\unskip\dimen5=\ht0 \advance\dimen5 \dimen
1 \dimen4=\dp0 \penalty2 \reformatrow\unpenalty\ht4=\dimen5 \dp4=\dimen4 \ht3%
\z@\dp3\z@\setbox1=\vbox{\box4 \nointerlineskip\box3 \nointerlineskip\unvbox1%
}\dimen2=\dimen1 \repeat\unvbox1}

\newif\ifcontinuerow

\def\reformatrow{\cdbfl=\cdbfr\penalty4 \noindent\unhbox0 \loopa\unskip
\setbox\cdbfl=\lastbox\ifhbox\cdbfl\advance\cdbfl1\repeat\endgraf\unskip
\unpenalty\dimen6=2\DiagramCellWidth\dimen7=-\DiagramCellWidth\setbox3=\hbox{%
}\setbox4=\hbox{}\setbox7=\box\voidb@x\cdvd=\cdvdl\continuerowtrue\loopa
\advance\cdvd-1 \adjustcells\ifcontinuerow\advance\dimen6\wd\cdbfl\cdda=.5%
\dimen6 \ifdim\cdda<\DiagramCellWidth\then\dimen6\DiagramCellWidth\advance
\dimen6-\cdda\nopendvert\cdda\DiagramCellWidth\fi\advance\dimen7\cdda\dimen6=%
\wd\cdbfl\reformatcell\advance\cdbfl-1 \repeat\advance\dimen7.5\dimen6
\outHarrow} \def\adjustcells{\ifnum\cdbfr>\cdbfl\then\ifnum\cdvdr>\cdvd\then
\continuerowfalse\else\setbox\cdbfl=\hbox to\wd\cdvd{\lcdtempl\VonH{}%
\rcdtempl}\fi\else\ifnum\cdvdr>\cdvd\then\advance\cdvdr-1 \setbox\cdvd=\vbox{%
}\wd\cdvd=\wd\cdbfl\dp\cdvd=\dp1 \fi\fi}

\def\reformatcell{\sethorizhtdp\penalty5 \noindent\unhbox\cdbfl\skip0=%
\lastskip\unskip\endgraf\ifcase\prevgraf\reformatempty\or\reformatobject\else
\reformatcomplex\fi\unskip\unpenalty}\def\reformatobject{\setbox6=\lastbox
\unskip\vadjdon6\outVarrow\setbox6=\hbox{\unhbox6}\advance\dimen7-.5\wd6
\outHarrow\dimen7=-.5\wd6 \pusht6\to4}\newcount\globnum

\def\reformatcomplex{\setbox6=\lastbox\unskip\setbox9=\lastbox\unskip\setbox9%
=\hbox{\unhbox9 \skip0=\lastskip\unskip\global\globnum\lastpenalty\hskip\skip
0 }\advance\globnum9999 \ifcase\globnum\reformathoriz\or\reformatpile\or
\reformatHonV\or\reformatVonH\or\reformatvert\or\reformatHmeetV\fi}

\def\reformatempty{\vpassdon\ifdim\skip0>\z@\then\hpassdon\else\ifvoid2 \then
\else\advance\dimen7-.5\dimen0 \cdda=\wd2\advance\cdda.5\dimen0\wd2=\cdda\fi
\fi}\def\VonH{\edef\HmeetingV{\string\VonH}\doVonH6}\def\HonV{\edef\HmeetingV
{\string\HonV}\doVonH7}\def\HmeetV{\edef\HmeetingV{\string\HmeetV}\MapBreadth
-2\MapShortFall\doVonH4}\def\doVonH#1{\cdna-999#1 \def\next{\mapctxterr
\HmeetingV}\ifallowhorizmap\ifallowvertmap\then\def\next{\futurelet
\thenexttoken\dooVonH}\allowhorizmapfalse\allowvertmapfalse\edef\whyforbidmap
{after \HmeetingV\space}\fi\fi\next}\def\dooVonH{\sethorizhtdp\ifx
\thenexttoken[\then\let\next\VonHstrut\else\sethorizhtdp\dimen0\MapBreadth
\let\next\VonHnostrut\fi\next}\def\VonHstrut[#1]{\setbox0=\hbox{$#1$}\dimen0%
\wd0\dimen8\ht0\dimen9\dp0 \VonHnostrut}\def\VonHnostrut{\setbox0=\hbox{}\ht0%
=\dimen8\dp0=\dimen9\wd0=.5\dimen0 \copy0\penalty\cdna\box0 }\def
\reformatHonV{\hpassdon\doreformatHonV}\def\reformatHmeetV{\dimen@=\wd9
\advance\dimen7-\wd9 \outHarrow\setbox6=\hbox{\unhbox6}\dimen7-\wd6 \advance
\dimen@\wd6 \setbox6=\hbox to\dimen@{\hss}\pusht6\to4\doreformatHonV}\def
\doreformatHonV{\setbox9=\hbox{\unhbox9 \unskip\unpenalty\global\setbox
\globbox=\lastbox}\vadjdon\globbox\outVarrow}\def\reformatVonH{\vpassdon
\advance\dimen7-\wd9 \outHarrow\setbox6=\hbox{\unhbox6}\dimen7=-\wd6 \setbox6%
=\hbox{\kern\wd9 \kern\wd6}\pusht6\to4}\def\hpassdon{}\def\vpassdon{\dimen@=%
\dp\cdvd\advance\dimen@\dimen4 \advance\dimen@\dimen5 \dp\cdvd=\dimen@
\nopendvert}\def\vadjdon#1{\dimen8=\ht#1 \dimen9=\dp#1 }

\def\HorizontalMap#1#2#3#4#5{\sethorizhtdp\setbox1=\makeharrowpart{#1}\def
\arrowfillera{#2}\def\arrowfillerb{#4}\setbox5=\makeharrowpart{#5}\ifx
\arrowfillera\justhorizline\then\def\arra{\hrule\horizhtdp}\def\kea{\kern-0.%
01em}\let\arrstruthtdp\horizhtdp\else\def\kea{\kern-0.15em}\setbox2=\hbox{%
\kea${\arrowfillera}$\kea}\def\arra{\copy2}\def\arrstruthtdp{height\ht2 depth%
\dp2 }\fi\ifx\arrowfillerb\justhorizline\then\def\arrb{\hrule\horizhtdp}\def
\keb{kern-0.01em}\ifx\arrowfillera\empty\then\let\arrstruthtdp\horizhtdp\fi
\else\def\keb{\kern-0.15em}\setbox4=\hbox{\keb${\arrowfillerb}$\keb}\def\arrb
{\copy4}\ifx\arrowfilera\empty\then\def\arrstruthtdp{height\ht4 depth\dp4 }%
\fi\fi\setbox3=\makeharrowpart{{#3}\vrule width\z@\arrstruthtdp}\ifincommdiag
\then\ifallowhorizmap\then\allowvertmapfalse\allowhorizmapfalse\def
\whyforbidmap{after horizontal map }\else\mapctxterr{horizontal}\fi\fi\let
\execmap\execHorizontalMap\gettwoargs}\def\justempty{\empty}\def
\makeharrowpart#1{\hbox{\mathsurround\z@\offinterlineskip\def\next{#1}\ifx
\next\justempty\else$\mkern-1.5mu{\next}\mkern-1.5mu$\fi}}\def\justhorizline{%
-}

\def\execHorizontalMap{\dimen0=\wd6 \ifdim\dimen0<\wd7\then\dimen0=\wd7\fi
\dimen3=\wd3 \ifdim\dimen0<2em\then\dimen0=2em\fi\skip2=.5\dimen0
\ifincommdiag plus 1fill\fi minus\z@\advance\skip2-.5\dimen3 \skip4=\skip2
\advance\skip2-\wd1 \advance\skip4-\wd5 \kern\MapShortFall\box1 \xleaders
\arra\hskip\skip2 \vbox{\lineskiplimit\maxdimen\lineskip.5ex \ifhbox6 \hbox to%
\dimen3 {\hss\box6\hss}\fi\vtop{\box3 \ifhbox7 \hbox to\dimen3 {\hss\box7\hss
}\fi}}\ifincommdiag\kern-.5\dimen3\penalty-9999\null\kern.5\dimen3\fi
\xleaders\arrb\hskip\skip4 \box5 \kern\MapShortFall}

\def\reformathoriz{\vadjdon6\outVarrow\ifvoid7\else\mapclasherr\fi\setbox2=%
\box9 \wd2=\dimen7 \dimen7=\z@\setbox7=\box6 }

\def\resetharrowpart#1#2{\ifvoid#1\then\ifdim#2=\z@\else\setbox4=\hbox{%
\unhbox4\kern#2}\fi\else\ifhbox#1\then\setbox#1=\hbox to#2{\unhbox#1}\else
\widenpile#1\fi\pusht#1\to4\fi}\def\outHarrow{\resetharrowpart2{\wd2}\pusht2%
\to4\resetharrowpart7{\dimen7}\pusht7\to4\dimen7=\z@}

\def\pile#1{{\incommdiagtrue\let\pile\innerpile\allowvertmapfalse
\allowhorizmaptrue\def\whyforbidmap{inside pile }\baselineskip.5\PileSpacing
\lineskip\z@\lineskiplimit\z@\mathsurround\z@\tabskip\z@\let\\\pilecr\vcenter
{\halign{\hfil$##$\hfil\cr#1 \crcr}}}\ifincommdiag\then\ifallowhorizmap\then
\penalty-9998 \allowvertmapfalse\allowhorizmapfalse\def\whyforbidmap{after
pile }\else\mapctxterr{pile}\fi\fi}\def\pilecr{\cr}\def\innerpile#1{\noalign{%
\halign{\hfil$##$\hfil\cr#1 \crcr}}}

\def\reformatpile{\vadjdon9\outVarrow\ifvoid7\else\mapclasherr\fi\penalty3
\setbox9=\hbox{\unhbox9 \unskip\unpenalty\setbox9=\lastbox\unhbox9 \global
\setbox\globbox=\lastbox}\unvbox\globbox\setbox9=\vbox{}\setbox7=\vbox{}%
\loopb\setbox6=\lastbox\ifhbox6 \skip3=\lastskip\unskip\splitpilerow\repeat
\unpenalty\setbox9=\hbox{$\vcenter{\unvbox9}$}\setbox2=\box9 \dimen7=\z@}\def
\pilestrut{\vrule height\dimen0 depth\dimen3 width\z@}\def\splitpilerow{%
\dimen0=\ht6 \dimen3=\dp6 \penalty6 \noindent\unhbox6\unskip\setbox6=\lastbox
\unskip\unhbox6\endgraf\setbox6=\lastbox\unskip\ifcase\prevgraf\or\setbox6=%
\hbox\tozpt{\hss\unhbox6\hss}\ht6=\dimen0 \dp6=\dimen3 \setbox9=\vbox{\vskip
\skip3 \hbox to\dimen7{\hfil\box6}\nointerlineskip\unvbox9}\setbox7=\vbox{%
\vskip\skip3 \hbox{\pilestrut\hfil}\nointerlineskip\unvbox7}\or\setbox7=\vbox
{\vskip\skip3 \hbox{\pilestrut\unhbox6}\nointerlineskip\unvbox7}\setbox6=%
\lastbox\unskip\setbox9=\vbox{\vskip\skip3 \hbox to\dimen7{\pilestrut\unhbox6%
}\nointerlineskip\unvbox9}\fi\unskip\unpenalty}

\def\widenpile#1{\setbox#1=\hbox{$\vcenter{\unvbox#1 \setbox8=\vbox{}\loopb
\setbox9=\lastbox\ifhbox9 \skip3=\lastskip\unskip\setbox8=\vbox{\vskip\skip3
\hbox to\dimen7{\unhbox9}\nointerlineskip\unvbox8}\repeat\unvbox8 }$}}

\def\VerticalMap#1#2#3#4#5{\setbox1=\makevarrowpart{#1}\def\arrowfillera{#2}%
\setbox3=\makevarrowpart{#3}\def\arrowfillerb{#4}\setbox5=\makevarrowpart{#5}%
\ifx\arrowfillera\justverticalline\then\def\arra{\vrule width\MapBreadth}\def
\kea{\kern-0.05ex}\else\def\kea{\kern-0.35ex}\setbox2=\vbox{\kea
\makevarrowpart\arrowfillera\kea}\def\arra{\copy2}\fi\ifx\arrowfillerb
\justverticalline\then\def\arrb{\vrule width\MapBreadth}\def\keb{\kern-0.05ex%
}\else\def\keb{\kern-0.35ex}\setbox4=\vbox{\keb\makevarrowpart\arrowfillerb
\keb}\def\arrb{\copy4}\fi\ifallowvertmap\else\mapctxterr{vertical}\fi\let
\execmap\execVerticalMap\allowhorizmapfalse\def\whyforbidmap{after vertical }%
\gettwoargs}\def\justverticalline{|}\def\makevarrowpart#1{\hbox to\MapBreadth
{\offinterlineskip\hss$\kern\MapBreadth{#1}$\hss}}

\def\execVerticalMap{\setbox3=\makevarrowpart{\box3}\setbox0=\hbox{}\ht0=\ht3
\dp0\z@\ht3\z@\box6 \setbox8=\vtop spread2ex{\offinterlineskip\box3 \xleaders
\arrb\vfill\box5 \kern\MapShortFall}\dp8=\z@\box8 \kern-\MapBreadth\setbox8=%
\vbox spread2ex{\offinterlineskip\kern\MapShortFall\box1 \xleaders\arra\vfill
\box0}\ht8=\z@\box8 \ifincommdiag\then\kern-.5\MapBreadth\penalty-9995 \null
\kern.5\MapBreadth\fi\box7\hfil}

\newcount\colno\newdimen\cdda\newbox\globbox\def\reformatvert{\setbox6=\hbox{%
\unhbox6}\cdda=\wd6 \dimen3=\dp\cdvd\advance\dimen3\dimen4 \setbox\cdvd=\hbox
{}\colno=\prevgraf\advance\colno-2 \loopb\setbox9=\hbox{\unhbox9 \unskip
\unpenalty\dimen7=\lastkern\unkern\global\setbox\globbox=\lastbox\advance
\dimen7\wd\globbox\advance\dimen7\lastkern\unkern\setbox9=\lastbox\vtop to%
\dimen3{\unvbox9}\kern\dimen7 }\ifnum\colno>0 \ifdim\wd9<\PileSpacing\then
\setbox9=\hbox to\PileSpacing{\unhbox9}\fi\dimen0=\wd9 \advance\dimen0-\wd
\globbox\setbox\cdvd=\hbox{\kern\dimen0 \box\globbox\unhbox\cdvd}\pushh9\to6%
\advance\colno-1 \setbox9=\lastbox\unskip\repeat\advance\dimen7-.5\wd6
\advance\dimen7.5\cdda\advance\dimen7-\wd9 \outHarrow\dimen7=-.5\wd6 \advance
\dimen7-.5\cdda\pusht9\to4\pusht6\to4\nopendvert\dimen@=\dimen6\advance
\dimen@-\wd\cdvd\advance\dimen@-\wd\globbox\divide\dimen@2 \setbox\cdvd=\hbox
{\kern\dimen@\box\globbox\unhbox\cdvd\kern\dimen@}\dimen8=\dp\cdvd\advance
\dimen8\dimen5 \dp\cdvd=\dimen8 \ht\cdvd=\z@}

\def\outVarrow{\ifhbox\cdvd\then\deepenbox\cdvd\pusht\cdvd\to3\else
\nopendvert\fi\dimen3=\dimen5 \advance\dimen3-\dimen8 \setbox\cdvd=\vbox{%
\vfil}\dp\cdvd=\dimen3} \def\nopendvert{\setbox3=\hbox{\unhbox3\kern\dimen6}}%
\def\deepenbox\cdvd{\setbox\cdvd=\hbox{\dimen3=\dimen4 \advance\dimen3-\dimen
9 \setbox6=\hbox{}\ht6=\dimen3 \dp6=-\dimen3 \dimen0=\dp\cdvd\advance\dimen0%
\dimen3 \unhbox\cdvd\dimen3=\lastkern\unkern\setbox8=\hbox{\kern\dimen3}%
\loopb\setbox9=\lastbox\ifvbox9 \setbox9=\vtop to\dimen0{\copy6
\nointerlineskip\unvbox9 }\dimen3=\lastkern\unkern\setbox8=\hbox{\kern\dimen3%
\box9\unhbox8}\repeat\unhbox8 }}

\newif\ifPositiveGradient\PositiveGradienttrue\newif\ifClimbing\Climbingtrue
\newcount\DiagonalChoice\DiagonalChoice1 \newcount\lineno\newcount\rowno
\newcount\charno\def\laf{\afterassignment\xlaf\charno='}\def\xlaf{\hbox{%
\tenln\char\charno}}\def\lah{\afterassignment\xlah\charno='}\def\xlah{\hbox{%
\tenln\char\charno}}\def\makedarrowpart#1{\hbox{\mathsurround\z@${#1}$}}\def
\lad{\afterassignment\xlad\charno='}\def\xlad{\hbox{\setbox2=\xlaf\setbox0=%
\hbox to.3\wd2{\hss.\hss}\dimen0=\ht0 \advance\dimen0-\dp0 \dimen1=.3\ht2
\advance\dimen1-\dimen0 \dp0=.5\dimen1 \dimen1=.3\ht2 \advance\dimen1\dimen0
\ht0=.5\dimen1 \raise\dp0\box0}}

\def\DiagonalMap#1#2#3#4#5{\ifPositiveGradient\then\let\mv\raise\else\let\mv
\lower\fi\setbox2=\makedarrowpart{#2}\setbox1=\makedarrowpart{#1}\setbox4=%
\makedarrowpart{#4}\setbox5=\makedarrowpart{#5}\setbox3=\makedarrowpart{#3}%
\let\execmap\execDiagonalLine\gettwoargs}

\def\makeline#1(#2,#3;#4){\hbox{\dimen1=#2\relax\dimen2=#3\relax\dimen5=#4%
\relax\vrule height\dimen5 depth\z@ width\z@\setbox8=#1\cdna=\dimen5 \divide
\cdna\dimen2 \ifnum\cdna=0 \then\box8 \else\dimen4=\dimen5 \advance\dimen4-%
\dimen2 \divide\dimen4\cdna\dimen3=\dimen1 \cdnb=\dimen2 \divide\cdnb1000
\divide\dimen3\cdnb\cdnb=\dimen4 \divide\cdnb1000 \multiply\dimen3\cdnb\dimen
6\dimen1 \advance\dimen6-\dimen3 \cdnb=0 \ifPositiveGradient\then\dimen7\z@
\else\dimen7\cdna\dimen4 \multiply\dimen4-1 \fi\loop\raise\dimen7\copy8 \ifnum
\cdnb<\cdna\hskip-\dimen6 \advance\cdnb1 \advance\dimen7\dimen4 \repeat\fi}}%
\newdimen\objectheight\objectheight1.5ex

\def\execDiagonalLine{\setbox0=\hbox\tozpt{\cdna=\xcoord\cdnb=\ycoord\dimen8=%
\wd2 \dimen9=\ht2 \dimen0=\cdnb\DiagramCellHeight\advance\dimen0-2%
\MapShortFall\advance\dimen0-\objectheight\setbox2=\makeline\box2(\dimen8,%
\dimen9;.5\dimen0)\setbox4=\makeline\box4(\dimen8,\dimen9;.5\dimen0)\dimen0=2%
\wd2 \advance\dimen0-\cdna\DiagramCellWidth\advance\dimen0 2\DiagramCellWidth
\dimen2\DiagramCellHeight\advance\dimen2-\MapShortFall\dimen1\dimen2 \advance
\dimen1-\ht1 \advance\dimen2-\ht2 \dimen6=\dimen2 \advance\dimen6.25\dimen8
\dimen3\dimen2 \advance\dimen3-\ht3 \dimen4=\dimen2 \dimen7=\dimen2 \advance
\dimen4-\ht4 \advance\dimen7-\ht7 \advance\dimen7-.25\dimen8
\ifPositiveGradient\then\hss\raise\dimen4\hbox{\rlap{\box5}\box4}\llap{\raise
\dimen6\box6\kern.25\dimen9}\else\kern-.5\dimen0 \rlap{\raise\dimen1\box1}%
\raise\dimen2\box2 \llap{\raise\dimen7\box7\kern.25\dimen9}\fi\raise\dimen3%
\hbox\tozpt{\hss\box3\hss}\ifPositiveGradient\then\rlap{\kern.25\dimen9\raise
\dimen7\box7}\raise\dimen2\box2\llap{\raise\dimen1\box1}\kern-.5\dimen0 \else
\rlap{\kern.25\dimen9\raise\dimen6\box6}\raise\dimen4\hbox{\box4\llap{\box5}}%
\hss\fi}\ht0\z@\dp0\z@\box0}

\newif\ifmoremapargs\def\gettwoargs{\setbox7=\box\voidb@x\setbox6=\box
\voidb@x\moremapargstrue\def\whichlabel{6}\def\xcoord{2}\def\ycoord{2}\def
\contgetarg{\def\whichlabel{7}\ifmoremapargs\then\let\next\getanarg\let
\contgetarg\execmap\else\let\next\execmap\fi\next}\getanarg}\def\getanarg{%
\futurelet\thenexttoken\switcharg}\def\getlabel#1#2#3{\setbox#1=\hbox{\kern.2%
em$\labelstyle{#3}$\kern.2em}\dimen0=\ht#1\advance\dimen0 .4ex\ht#1=\dimen0
\dimen0=\dp#1\advance\dimen0 .4ex\dp#1=\dimen0 \contgetarg}\def
\eatspacerepeat{\afterassignment\getanarg\let\junk= }\def\catcase#1:{{\ifcat
\noexpand\thenexttoken#1\then\global\let\xcase\docase\fi}\xcase}\def\tokcase#%
1:{{\ifx\thenexttoken#1\then\global\let\xcase\docase\fi}\xcase}\def\default:{%
\docase}\def\docase#1\esac#2\esacs{#1}\def\skipcase#1\esac{}\def
\getcoordsrepeat(#1,#2){\def\xcoord{#1}\def\ycoord{#2}\getanarg}\let\esacs
\relax\def\switcharg{\global\let\xcase\skipcase\catcase{&}:\moremapargsfalse
\contgetarg\esac\catcase\bgroup:\getlabel\whichlabel-\esac\catcase^:\getlabel
6\esac\catcase_:\getlabel7\esac\tokcase{~}:\getlabel3\esac\tokcase(:%
\getcoordsrepeat\esac\catcase{ }:\eatspacerepeat\esac\default:%
\moremapargsfalse\contgetarg\esac\esacs}

\catcode`\@=\cdna

\newcommand{\Z}{{\bf Z}}

\newcommand{\Q}{{\bf Q}}
\newcommand{\GQ}{G_{\bf Q}}
\newcommand{\HH}{{\rm H}}
\newcommand{\Hom}{{\rm Hom}}
\newcommand{\Pp}{{\bf P}}
\newcommand{\Gal}{{\rm Gal}}
\newcommand{\D}{{\mathfrak{D}}}

\newtheorem{theo}{Th\'eor\`eme}
\newtheorem{cor}{Corollaire}

\def\refname{Bibliographie}
\def\abstractname{R\'esum\'e}
\def\keywordsname{Mots clefs}

\begin{document}

\title[Familles de polyn\^omes et courbes elliptiques
quotient]{Familles de polyn\^omes li\'ees aux courbes modulaires
  $X_1(l)$ unicursales et points rationnels non-triviaux de courbes
  elliptiques quotient}

\author{\sc Franck Lepr\'evost}
\address{Universit\'e Joseph Fourier, UFR de Math\'ematiques, 100, rue des Maths - B.P. 74 - F-38402 St-Martin d'H\`eres Cedex, France}
\email{Franck.Leprevost@ujf-grenoble.fr}

\author{\sc Michael Pohst}
\address{Technische Universit\"at Berlin, Fakult\"at II - Mathematik
  MA 8-1 - Stra\ss{}e des 17. Juni 136, D-10623 Berlin, Allemagne} 
\email{pohst@math.tu-berlin.de}

\author{\sc Andreas Sch\"opp}
\address{Technische Universit\"at Berlin, Fakult\"at II - Mathematik
  MA 8-1 - Stra\ss{}e des 17. Juni 136, D-10623 Berlin, Allemagne}
\email{schoepp@math.tu-berlin.de}

\keywords{Courbes modulaires, extensions di\'edrales, extensions
  cycliques, unit\'es 
  \vspace*{1cm}\hspace*{7.5cm} \it Article submitted to Acta Arithmetica}


\begin{abstract}
Soit $l$ un entier et $E_{c,l}$ la famille de Kubert des courbes
elliptiques d\'efinies sur $\Q$ munies d'un point rationnel $A$
d'ordre $l$. On note $F_{c,l}$ la courbe elliptique quotient de
$E_{c,l}$ par le groupe engendr\'e par $A$, et
$\varphi_l$ l'isog\'enie de $E_{c,l}$ sur $F_{c,l}$. Pour $l=3,4,5$ et
$6$, nous construisons explicitement, pour des param\'etrisations
convenables de $c$, des \'el\'ements non-triviaux de
$F_{c,l}(\Q)/\varphi_{l}(E_{c,l}(\Q))$, autrement dit, des points
explicites de $F_{c,l}(\Q)$ qui ne sont l'image par $\varphi_{l}$
d'aucun \'el\'ement de $E_{c,l}(\Q)$. Ces points sont en g\'en\'eral
d'ordre infini. Nous donnons des applications de cette m\'ethode \`a la
construction d'extensions cycliques de $\Q$ de degr\'e $l$, et
retrouvons certains corps obtenus par Shanks et Gras. Dans un
article ult\'erieur, nous \'etudierons les propri\'et\'es
arithm\'etiques de certaines des extensions obtenues ici.
\end{abstract}

\maketitle

\section{Introduction}

L'\'etude du rang des courbes elliptiques d\'efinies sur $\Q$ et
${\Q}(t)$ a connu
ces derni\`eres ann\'ees de nombreux progr\`es, initi\'es par les
travaux de Mestre. En fait, il est (verbalement) conjectur\'e
l'existence de courbes elliptiques sur $\Q$ de grand rang et de groupe
de torsion rationnel arbitraire parmi la liste des quinze groupes
possibles, \`a savoir $\Z/l\Z$ pour $1 \leq l \leq 10$ ou $l=12$, ou
bien $\Z/2\Z \oplus \Z/2l\Z$ pour $1 \leq l \leq 4$. Dans plusieurs
travaux (\cite{fermigier}, \cite{kihara0}, \cite{kihara1},
\cite{odile}), l'on construit des courbes elliptiques de rang
\'elev\'e ayant un groupe de torsion non-r\'eduit \`a l'\'el\'ement
neutre. \\ 

\noindent Consid\'erons ici $3\leq l \leq 10$ ou $l=12$, et soit $E_{c,l}$ la
famille de Kubert (\cite{kubert}) des courbes elliptiques d\'efinies
sur $\Q$ munies 
d'un point rationnel $A$ d'ordre $l$ (dans le cas $l=3$, $c$ d\'esigne
en r\'ealit\'e un couple de param\`etres). En d'autres termes, $(E_{c,l},
A)$ param\'etrise les points rationnels sur $\Q$ de la courbe
modulaire $X_1(l)$, qui est isomorphe \`a ${\Pp}^1$ pour ces valeurs de
$l$ (elle est encore unicursale pour $l=1, 2$, valeurs qui n'ont
gu\`ere d'int\'er\^et dans le cadre qui nous occupe). Avec ces
notations, les travaux cit\'es plus haut exhibent, via des
param\'etrisations ing\'enieuses de $c$, des familles explicites de
courbes du type $E_{c,l}$ de $\Q$-rang diff\'erent de $0$. Les
m\'ethodes employ\'ees exploitent essentiellement le fait que le
param\`etre $c$ intervient avec un degr\'e relativement {\it petit}
dans des \'equations bien choisies de $E_{c,l}$. \\

\noindent Si $<A>$ d\'esigne le groupe engendr\'e par $A$, notons
$F_{c,l}$ la courbe elliptique quotient $E_{c,l}/<A>$, et 
$\varphi_l$ l'isog\'enie de $E_{c,l}$ sur $F_{c,l}$. On a ainsi la
suite exacte suivante de $\GQ$-modules galoisiens (o\`u $\GQ = \Gal
(\overline{\Q}/\Q)$) :
$$\begin{diagram}
1 & \rArr & \Z/l\Z \simeq <A> & \rArr & E_{c,l}({\overline{\Q}}) &
\rArr^{\varphi_{l}} & F_{c,l}({\overline{\Q}}) & \rArr & 1. \\
\end{diagram}$$

\noindent Connaissant une \'equation de $E_{c,l}$ et les coordonn\'ees
de $A$, les formules de V\'elu (\cite{velu}) permettent d'obtenir
explicitement une \'equation de $F_{c,l}$ et de $\varphi_l$. \\

\noindent Bien \'evidemment, l'image par $\varphi_l$ de
points de $E_{c,l}(\Q)$ fournit des points de $F_{c,l}(\Q)$, que nous
appelons ici des points rationnels triviaux de la courbe elliptique
quotient $F_{c,l}$. \\

\noindent Le probl\`eme, auquel nous nous int\'eressons dans la partie
$2$, consiste en la construction, via des param\'etrisations
convenables ou des sp\'ecialisations de $c$, d'\'el\'ements
non-triviaux de $F_{c,l}(\Q)/\varphi_{l}(E_{c,l}(\Q))$, autrement dit,
de points explicites de $F_{c,l}(\Q)$ qui ne sont l'image par
$\varphi_l$ d'aucun \'el\'ement de $E_{c,l}(\Q)$. En g\'en\'eral les
repr\'esentants de ces points dans $F_{c,l}(\Q)$ sont d'ordre
infini. Ceci dit, et bien que les courbes elliptiques $E_{c,l}(\Q)$ et
$F_{c,l}(\Q)$ aient m\^eme rang, nous n'utilisons pas les
constructions de \cite{fermigier}, \cite{kihara0}, \cite{kihara1},
\cite{odile} qui ne fourniraient, dans notre terminologie, que des points
triviaux de $F_{c,l}(\Q)$, comme remarqu\'e plus haut. Le probl\`eme
trait\'e ici devient rapidement d\'elicat (au regard de $l$), puisque
le degr\'e du 
param\`etre $c$ dans les \'equations de $F_{c,l}$ est, comme nous le
verrons en particulier dans les cas consid\'er\'es, notablement plus
\'elev\'e que dans les \'equations de $E_{c,l}$. Dans la partie $3$, nous
construisons, \`a partir de la donn\'ee $(E_{c,l},A)$, un polyn\^ome
$P_{n,c,l} \in {\Z}[n,c][x]$, pour lequel nous montrons que son corps
de d\'ecomposition sur ${\Q}(n,c)$ est g\'en\'eriquement le groupe
di\'edral $D_l$ \`a $2l$ \'el\'ements. Dans la partie $4$, nous
consid\'erons plus sp\'ecifiquement le cas $l=5$, et montrons que
notre construction permet de retrouver la famille g\'en\'erique de
Brumer, qui est isomorphe \`a celle obtenue, ind\'ependamment, par
Darmon (voir les r\'ef\'erences donn\'ees dans cette partie). Il est
tentant de regarder sous quelles conditions sur les param\`etres
$(n,c)$ le groupe de Galois de $P_{n,c,l}$ sur ${\Q}(n,c)$ devient
isomorphe au groupe ${\Z}/l{\Z}$. Dans la partie $5$, nous montrons
que ces conditions reviennent \`a la construction d'\'el\'ements
non-triviaux de $F_{c,l}(\Q)/\varphi_{l}(E_{c,l}(\Q))$. Les
r\'esultats de la partie $2$ permettent alors de construire
explicitement de telles extensions cycliques. Nous retrouvons
\'egalement les {\it simplest cubic fields} de Shanks et les {\it simplest
  quartic fields} de M.-N. Gras. Dans \cite{fam2}, nous
\'etudierons des propri\'et\'es arithm\'etiques d'extensions
construites ici. \\ 

\noindent Les calculs effectu\'es pour cet article ont n\'ecessit\'e
un usage tr\`es important des logiciels de calcul formel KANT
(\cite{kant}), MAGMA (\cite{magma}), MAPLE (\cite{maple}) et PARI
(\cite{pari}).

\section{Points non-triviaux de courbes elliptiques quotient}

\noindent Nous montrons ici le r\'esultat suivant :  

\begin{theo}
Pour $l=3,4,5$ et $6$, nous construisons explicitement, pour des
param\'etrisations convenables de $c$, des \'el\'ements non-triviaux de
$F_{c,l}(\Q)/\varphi_{l}(E_{c,l}(\Q))$. Ces points sont en g\'en\'eral
d'ordre infini.
\end{theo}

\noindent Par un argument de cohomologie galoisienne classique, comme
$\GQ = \Gal (\overline{\Q}/\Q)$ op\`ere trivialement sur $<A>$, on
\'etablit ais\'ement le : 
\begin{cor}
Pour $l=3, 4, 5$ et $6$, et pour les param\'etrisations de $c$ du
th\'eor\`eme pr\'ec\'edent, le groupe de cohomologie ${\HH}^1(\GQ,
\Z/l\Z) = {\Hom}(\GQ,\Z/l\Z)$ est non
r\'eduit \`a $0$. 
\end{cor}

\noindent Bien entendu, la suite exacte longue en cohomologie se
poursuivant, nous obtenons de la sorte plus pr\'ecis\'ement des
\'el\'ements non-triviaux du noyau de l'application 

$$\begin{diagram}
{\HH}^1(\GQ, \Z/l\Z) = {\Hom}(\GQ,\Z/l\Z) & \rArr & {\HH}^1(\GQ,
E_{c,l}({\overline{\Q}})). \\  
\end{diagram}$$

\subsection{Les cas $l=3$ et $l=5$}

\noindent Etant donn\'e que notre strat\'egie est suffisament
g\'en\'erale, nous la d\'eveloppons pour $l \geq 5$ impair, la
d\'etaillons sur le cas $l=5$, et serons plus
succints pour le cas $l=3$. Nous indiquons dans une remarque les
limites de notre approche. La diff\'erence de traitement du cas $l=3$
correspond essentiellement au fait que la famille des courbes
elliptiques d\'efinies sur ${\bf Q}$ munies d'un point d'ordre $3$ est
param\'etris\'ee par deux param\`etres, et non pas un, ce qui alourdit
quelque peu les notations. \\

\noindent Soit donc $l \geq 5$ impair. Les formules de
V\'elu permettent d'obtenir une \'equation de $F_{c,l}$ de la forme : 
$$y^2 = f_{c,l}(x) = 4x^3 + {\alpha}_l(c)x^2 + {\beta}_l(c) x +
{\gamma}_l(c),$$ 
o\`u $\alpha_l, \beta_l$ et $\gamma_l$ sont les \'el\'ements de
$\Z[c]$ donn\'es dans le tableau suivant pour $l=5$ : 

\medskip
\begin{center}
 \begin{tabular}{|c|c|} \hline
   $\alpha_5(c)$ & $c^2-30c+1$\\ \hline
   $\beta_5(c)$  & $-2c(3c+1)(4c-7)$ \\ \hline
   $\gamma_5(c)$ & $-c(4c^4-4c^3-40c^2+91c-4)$ \\ \hline
 \end{tabular}
\end{center}
\medskip

\noindent On cherche alors $x_{c,l}$ sous la forme d'un
polyn\^ome en $c$, de petit degr\'e, de sorte que l'on ait une identit\'e :
$$f_{c,l}(x_{c,l}) = A_l(c)G_l^2(c),$$
o\`u $A_l, G_l$ sont des polyn\^omes en $c$ tels que le degr\'e de $A_l$ en
$c$ soit \'egal \`a $1$ ou \`a $2$. Pour r\'ealiser cela, on calcule le
discriminant par rapport \`a $c$ de $f_{c,l}(x_{c,l})$, qui est un
polyn\^ome en les coefficients de $c$ de $x_{c,l}$. Il s'agit alors de
l'annuler. De tels $x_{c,l}$, et les $G_l(c)$ et $A_l(c)$ associ\'es, sont
donn\'es dans le tableau suivant pour $l=5$ :  

\medskip
\begin{center}
 \begin{tabular}{|c|c|} \hline
   $x_{c,5}$   & $-(u_0+1)c^2+(11u_0+8)c+u_0$ \\ \hline 
   $G_{5}(c)$  & $c^2 - 11c - 1$ \\ \hline
   $A_5(c)$    & $-(4u_0+3)(u_0+1)^2c^2 + 2(2u_0+1)(11u_0^2+11u_0+2)c
   + u_0^2(4u_0+1)$ \\ \hline
 \end{tabular}
\end{center}
\medskip

\noindent Il reste alors \`a param\'etriser, \`a l'aide d'un choix
appropri\'e des param\`etres restants la conique en $(c,z)$
d'\'equation :
$$A_l(c) = z^2.$$

\noindent {\underline{Cas $l = 5$}} : Comme nous l'avons soulign\'e
plus haut, le degr\'e en $c$ de $A_5(c)$ est $\leq 2$. Par
cons\'equent, plusieurs cas se pr\'esentent, selon que l'on choisisse
le degr\'e en $c$ de $A_5(c)$ \'egal \`a $1$ ou \`a $2$, et qui
peuvent simplifier les calculs. Le tableau suivant r\'esume les choix
et calculs effectu\'es. Les deux premi\`eres lignes correspondent \`a
l'annulation du coefficient de degr\'e $2$ en $c$ de $A_5(c)$, et la
troisi\`eme assure que le point de coordonn\'ees $(c,z) = (0,u_{0}t)$
appartient \`a la conique d'\'equation $A_5(c) = z^2$. 

\medskip
\begin{center}
 \begin{tabular}{|c|c|c|} \hline
  $u_0 = -1$ & $c = \frac{z^2-3}{4}$ & $x_{c,5} = \frac{5-3z^2}{4}$ \\ \hline 
  $u_0 = -\frac{3}{4}$ & $c = 16z^2+18$ & $x_{c,5} = -64z^4-148z^2-\frac{345}{4}$ \\ \hline 
  $u_0 = \frac{t^2-1}{4}$ & $c =
  \frac{11t^6+33t^4-8mt^3+21t^2+8mt-1}{t^6+8t^4+21t^2+16m^2+18}$ &
  $x_{c,5} = -\frac{t^2+3}{4} c^2 + \frac{11t^2+21}{4} c +
\frac{t^2-1}{4}$ \\ \hline 
 \end{tabular}
\end{center}
\medskip

\noindent Il d\'ecoule de ce qui pr\'ec\`ede que, pour ces choix des
param\`etres, la courbe elliptique quotient $F_{c,l}(\Q)$ poss\`ede le
point de coordonn\'ees $(x_{c,l}, zG_l(c))$. Le calcul montre
que le point ainsi construit est g\'en\'eriquement d'ordre infini. 
Le calcul montre \'egalement que ce
point n'est g\'en\'eriquement pas l'image par l'isog\'enie $\varphi_l$
d'un point rationnel de $E_{c,l}$. Cette derni\`ere question sera
consid\'er\'ee de nouveau dans la partie $5$. \\

\noindent {\underline{Cas $l = 3$}} : Les formules de
V\'elu permettent d'obtenir une \'equation de $F_{a_1,a_3,l}$ de la forme : 
$$y^2 = f_{a_1,a_3,3}(x) = 4x^3 + a_1^2x^2 -18a_1a_3x - a_3(4a_1^3+27a_3).$$ 
Si l'on choisit $x_{a_1,a_3,3} = u_1a_3+\frac{a_1u_1+1}{u_1^2},$ le
calcul montre que 
$$f_{a_1,a_3,3}(x_{a_1,a_3,3}) = A_{3}(a_1,a_3)G^2_{3}(a_1,a_3),$$
o\`u $A_{3}(a_1,a_3) = 4u_1^3a_3+(u_1a_1+1)^2$ et $G_{3}(a_1,a_3) =
\frac{u_1^3a_3-u_1a_1-2}{u_1^3}$. L'\'equation $A_l(c) = z^2$ est
lin\'eaire en $a_3$, et il suffit de prendre $a_3 =
\frac{z^2-(u_1a_1+1)^2}{4u_1^3}$. On obtient alors $x_{a_1,a_3,3} =
\frac{z^2-(a_1u_1+1)(a_1u_1-3)}{4u_1^2}$, et l'on conclut comme dans
le cas $l=5$.\\

\noindent Remarque : Pour les autres valeurs impaires de $l$, la m\'ethode
d\'ecrite ici trouve ses limites essentiellement dans le calcul du
discriminant de $f_{c,l}(x_{c,l})$ par rapport \`a $c$. Les logiciels
de calcul formel ne permettent pas de factoriser en toute
g\'en\'eralit\'e cette quantit\'e. Cependant, on constate que celle-ci
s'exprime comme un produit d'un {\it gros} facteur par un {\it petit}
facteur, ce dernier intervenant \`a la puissance $l$. On peut, par
sp\'ecialisation et interpolation, calculer explicitement ce petit
facteur, du moins l'avons nous fait pour $l=7$. Nous ne sommes en
revanche malheureusement pas parvenus \`a l'annuler de mani\`ere utile
dans le cadre que nous consid\'erons ici. \\

\subsection{Les cas $l= 4$ et $l = 6$}

Dans le cas o\`u $l$ est pair, l'approche est l\'eg\`erement
diff\'erente. Cependant, l\`a \'egalement, nous la d\'ecrivons pour $l
\geq 6$ pair, la d\'etaillons sur le cas $l=6$, et serons plus
succints pour le cas $l=4$. \\

\noindent Les formules de V\'elu permettent de nouveau
d'obtenir une \'equation de $F_{c,l}$, qui est de la forme :
$$y^2 = f_{c,l}(x) = (4x - \alpha_l(c)) (x^2 + \beta_l(c)x + \gamma_l(c)),$$
o\`u $\alpha_l(c), \beta_l(c), \gamma_l(c)$ sont des \'el\'ements de
$\Z[c]$ donn\'es dans le tableau suivant pour $l=6$ :

\medskip
\begin{center}
 \begin{tabular}{|c|c|} \hline
   $\alpha_6(c)$ & $19c^2+14c-1$\\ \hline
   $\beta_6(c)$  & $2c(2c+1)$ \\ \hline
   $\gamma_6(c)$ & $c(4c^3+4c^2+c+4)$ \\ \hline
 \end{tabular}
\end{center}
\medskip

\noindent On cherche de nouveau $x_{c,l}$ sous la forme d'un
polyn\^ome en $c$ de sorte que l'on ait une identit\'e 
$$f_{c,l}(x_{c,l}) = A_l(c) G_l(c)^2,$$
o\`u $A_l$ est un polyn\^ome en $c$ de degr\'e $\leq 2$. Mais pour cela, on
exploite tout d'abord la factorisation de $f_{c,l}(x)$, en prenant 
$x_{c,l} = \frac{\alpha_l(c) + u^2}{4},$
ce qui assure que $4x_{c,l} - \alpha_l(c) = u^2.$ Ensuite, on cherche
$u$ sous la forme d'un polyn\^ome de degr\'e petit en $c$, de sorte que 
$x_{c,l}^2 + \beta_l(c)x_{c,l} + \gamma_l(c)$ admette de nouveau des
facteurs carr\'es. En pratique, on prend $u = v_1 c + v_0$, et l'on
cherche \`a sp\'ecialiser les param\`etres $v_0, v_1$ pour avoir la
factorisation $f_{c,l}(x_{c,l}) = A_l(c) G_l(c)^2$ voulue. De tels
$x_{c,l}$, et les $G_l(c)$ et $A_l(c)$ associ\'es, sont donn\'es dans
le tableau suivant pour $l=6$ :   

\medskip
\begin{center}
 \begin{tabular}{|c|c|} \hline
   $x_{c,6}$   & $\frac{19c^2+14c-1+v_0^2(9c+1)^2}{4}$ \\ \hline 
   $G_{6}(c)$  & $\frac{v_0(9c+1)^2}{4}$ \\ \hline
   $A_6(c)$    & $9(3v_0^2+1)^2c^2+2(3v_0^2+1)(3v_0^2+5)c+(v_0^2-1)^2$
   \\ \hline 
 \end{tabular}
\end{center}
\medskip

\noindent {\underline{Cas $l = 6$}} : La conique d'\'equation $A_6(c)
= z^2$ contient le point rationnel de coordonn\'ees $(c,z) = (0,
v_0^2-1)$, et donc se param\'etrise, et l'on trouve finalement 
$$c= 2 \frac{9v_0^4+18v_0^2-v_0^2z+z+5}{(z+3+9v_0^2)(z-3-9v_0^2)}.$$
On v\'erifie alors que le point de $F_{c,6}$ d'abscisse
$x_{c,6}$ ainsi construit ne provient pas d'un point de $E_{c,6}$ via
l'isog\'enie $\varphi_6$.\\

\noindent {\underline{Cas $l = 4$}} : Dans ce cas, les formules de
V\'elu donnent une \'equation de $F_{4,c}$ de la forme :
$$y^2 = f_{c,4}(x) = (x+c)(4x^2+x+c).$$
On choisit $x_{c,4} = u^2 - c,$ si bien que 
$$f_{c,4}(x_{c,4}) = u^2(4c^2-8u^2c+u^2(4u^2+1)).$$
La conique d'\'equation $4c^2-8u^2c+u^2(4u^2+1) = z^2$ se
param\'etrise ais\'ement, et l'on trouve 
$$c= \frac{u^2(4u^2+1)-v^2}{4v+8u^2}.$$
De m\^eme, on v\'erifie alors que le point de $F_{c,4}$ d'abscisse
$x_{c,4}$ ainsi construit ne provient pas d'un point de $E_{c,4}$ via
l'isog\'enie $\varphi_4$.

\section{Courbes elliptiques et polyn\^omes \`a groupe de Galois di\'edral}

\noindent Avec les notations de la partie pr\'ec\'edente, soit
$P_{n,c,l}(x)$ le polyn\^ome de $\Z[n,c][x]$ d\'efini par  
$$P_{n,c,l}(x) = \prod_{i=0}^{l-1}(x-x(P+iA)) = x^l - nx^{l-1} +
\cdots,$$ o\`u $P$ d\'esigne un point non $\Q$-rationnel de $E_{c,l}$,
$A$ un point fix\'e d'ordre $l$, et
$x(P+iA)$ l'abscisse du point $P+iA \in E_{c,l}$. Les logiciels de
calcul formel permettent d'obtenir l'\'equation explicite \footnote{On
  peut r\'ecup\'erer les \'equations de $P_{n,c,l}(x)$ sur \\
  http://www.math.tu-berlin.de/$\sim$
  kant/publications/papers/polynomes.txt} de $P_{n,c,l}(x)\in
\Z[n,c][x]$, et d'\'etablir le r\'esultat suivant :  

\begin{theo}
Soit $l$ un entier tel que $3\leq l \leq 10$ ou $l=12$. Le polyn\^ome
$P_{n,c,l}$ construit ci-dessus est g\'en\'eriquement irr\'eductible
sur le corps ${\Q}(n,c)$, et le groupe de Galois sur ${\Q}(n,c)$ de
son corps de d\'ecomposition est g\'en\'eriquement le groupe di\'edral
$D_l$ \`a $2l$ \'el\'ements. 
\end{theo}

\section{Le cas $D_5$}

\noindent Plusieurs auteurs se sont int\'eress\'es \`a la construction
de polyn\^omes quintiques de groupe de Galois $D_5$. Ainsi Weber
(\cite{weber}, p. 676) et Cebotarev (\cite{cebo}, p. 344) donnent une 
condition n\'ecessaire et suffisante, sous la forme d'une
param\'etrisation explicite des coefficients $a$ et $b$, pour que $x^5
+ ax + b$ soit r\'esoluble par radicaux. Une telle caract\'erisation,
apparament obtenue de mani\`ere ind\'ependante, est \'egalement
l'objet de l'article \cite{will}. A partir de la caract\'erisation
d\^ue \`a Weber et Cebotarev, Roland, Yui et Zagier (\cite{zagier})
donnent la param\'etrisation des polyn\^omes quintiques $x^5 + ax + b$
ayant $D_5$ pour groupe de Galois. D'autres auteurs se sont
int\'eress\'es \`a ces questions (sans pr\'etendre en aucune mani\`ere
\`a l'exclusivit\'e, citons \cite{yui} pour les groupes $D_p$, o\`u
$p$ est premier, et \cite{mike1} pour une th\'eorie sur les relations
entre tours modulaires et groupes di\'edraux). \\

\noindent Dans le cas particulier $l=5$, la construction d\'ecrite
dans la partie 3 donne le polyn\^ome  
$$P_{n,c,5}(x) =
x^5-nx^4-(-c^3-2nc+c^2+c)x^3-(c^3+nc^2-3c^2)x^2-(-c^4+3c^3)x+c^4.$$
La substitution $(x,n,c) \longrightarrow (\frac{s}{x}, -u, s)$ redonne
la famille g\'en\'erique de Brumer (\cite{armand}) cit\'ee par Martinais
et Schneps (\cite{masch}, p. 151) : 
$$B_{s,u}(x) = x^5 + (s-3) x^4 + (u-s+3)x^3 + (s^2-s-2u-1)x^2 + ux +
s.$$

\noindent Cette famille est g\'en\'erique dans le sens o\`u Brumer
 (\cite{masch}, p. 151) affirme que, si
$F$ est un corps contenant $\Q$ et $K$ est une extension galoisienne
de $F$ de groupe de Galois $D_5$, alors $K$ est le corps de
d\'ecomposition d'un polyn\^ome de la forme $B_{s,u}(x)$ pour des
valeurs de $s$ et $u$ appartenant \`a $F$. Malheureusement, \`a
l'heure actuelle, nous ne disposons pas de la preuve de ce fait. Par
ailleurs, Kihel (\cite{omar}, p. 471)
rappelle la construction de Darmon (\cite{darmon}) de la famille  
$$D_{S,T}(x) = x^5 -Sx^4 + (T+S+5)x^3 - (S^2+S-2T-5)x^2+(T+2S+5)x
-(S+3).$$
Cette famille est encore isomorphe \`a celle de Brumer, comme on le
constate \`a l'aide de la
transformation $(x,s,u) \longrightarrow (-x,S+3,T+2S+5)$. Les
constructions de Brumer, de Darmon et celle pr\'esent\'ee ici
produisent donc la m\^eme famille de polyn\^omes, obtenue de
mani\`ere ind\'ependante par les diff\'erents auteurs : nous avons
d\'ecouvert l'article \cite{omar} et l'existence des notes
\cite{darmon}, dont l'original ne semble plus disponible
\cite{darmon2} mais que l'article \cite{omar} d\'ecrit pour
l'essentiel, apr\`es avoir d\'emontr\'e le th\'eor\`eme 2 en toute
g\'en\'eralit\'e, ce qui inclut en particulier le cas $l=5$. Le cas
$l=5$ est \'egalement repris dans \cite{omar2}. Par
ailleurs, nous n'avons malheureusement pas eu d'informations
concr\`etes concernant la construction de Brumer et la preuve de son
r\'esultat de g\'en\'ericit\'e, qui n'est explicite ni dans
\cite{masch}, ni dans \cite{armand}. La question d'\'etendre le
r\'esultat de Brumer aux autres cas, c'est-\`a-dire de d\'ecrire
explicitement les extensions di\'edrales d'ordre $2l$ de $\Q$ reste
donc {\it a priori} encore ouverte pour les cas $l \neq 5$. 

\section{Une application : construction de certaines extensions
  cycliques de $\Q$} 

\noindent Pour $3 \leq l \leq 10$ ou $l=12$, il parait naturel de
chercher, par sp\'ecialisation des param\`etres dans $P_{n,c,l}(x)$,
des familles ou des exemples de polyn\^omes dont le groupe de Galois
est isomorphe \`a $\Z/l\Z$. Nous montrons ici les r\'esultats suivants
:  

\begin{theo}
Pour $3 \leq l \leq 6$, il existe une famille explicite de polyn\^omes
de degr\'e $l$, d\'efinie sur $\Q$, \`a groupe de Galois cyclique
d'ordre $l$. Ces familles sont index\'ees sur, d'une part, un
param\`etre rationnel, d'autre part, un point d'ordre infini d'une
famille de courbes elliptiques.
\end{theo}

\begin{theo}
Pour $7 \leq l \leq 10$ et $l=12$, il existe une famille explicite de
polyn\^omes de degr\'e $l$, d\'efinie sur $\Q$, \`a groupe de Galois
cyclique d'ordre $l$. Ces familles sont index\'ees par un point
d'ordre infini sur une courbe elliptique d\'efinie sur $\Q$. 
\end{theo}

\noindent La strat\'egie que nous adoptons est la
suivante : sachant que le groupe de Galois de $P_{n,c,l}$ est
g\'en\'eriquement $D_{l}$, il existe un \'el\'ement $\D_l$, appel\'e
fonction d'indicateur de l'extension, qui s'exprime a priori en
fonction des racines de $P_{n,c,l}$, et qui satisfait une equation
quadratique : 
$$\D_l^2 + U_l(n,c) \D_l + V_l(n,c) = 0.$$

\noindent Le groupe de Galois de $P_{n,c,l}$ est isomorphe \`a
$\Z/l\Z$ si et seulement si $\D_l$ est un rationnel, et si $P_{n,c,l}$
est irr\'eductible. La condition de rationnalit\'e de $\D_l$
\'equivaut pr\'ecis\'ement \`a trouver un \'el\'ement de
$F_{c,l}(\Q)$. La condition d'irr\'eductibilit\'e de $P_{n,c,l}$
\'equivaut \`a ce que cet \'el\'ement de $F_{c,l}(\Q)$ ne soit l'image
par $\varphi_{c,l}$ d'aucun \'el\'ement de $E_{c,l}(\Q)$. Les
r\'esultats de la partie $2$ permettent de conclure. \\

\noindent Cette approche est coh\'erente avec l'interpr\'etation
cohomologique donn\'ee dans la partie $2$. En effet, \`a un point
ferm\'e de ${\Hom}(\GQ, \Z/l\Z)$ correspond une extension galoisienne
de $\Q$ \`a groupe de Galois $\Z/l\Z$, du moins si $l$ est premier,
celle-ci s'obtenant comme la fibre de $\varphi_l$. \\

\noindent Il est \`a noter que plusieurs familles remarquables de
polyn\^omes existent fournissant des extensions cycliques de $\Q$. Les
corps cubiques 
les plus simples ont ainsi \'et\'e introduits par D. Shanks
(\cite{dan}). E. Lehmer (\cite{emma}) a construit des corps quintiques
cycliques simples, et M.-N. Gras des corps quartiques cycliques
simples (\cite{mng1}) et des corps sextiques cycliques simples
(\cite{mng2}). D'autres extensions cycliques de degr\'e $6$ et $10$
ont \'et\'e consid\'er\'ees par O. Lecacheux dans, respectivement,
\cite{odile1} et \cite{odile2}. \\ 

\noindent Notre construction permet de retrouver certaines de
ces extensions {\it simples}. En
effet, dans le cas $l=3$, le corps engendr\'e par $P_{n,c,3}(x)$ l'est
aussi (via des transformations \'el\'ementaires constitu\'ees de
translation et homoth\'etie) par
${\tilde{P}}_{n,c,3}(x)= x^3 +ux^2 -nx +v$. 
Il suffit alors de prendre $(u,v,n) = (-t,-1,t+3)$ pour retrouver la
famille cubique cyclique 
$$X^3-tX^2-(t+3)X-1$$
de Shanks. Dans le cas $l=4$, on constate, \`a l'aide de
transformations \'el\'ementaires (translation et homoth\'etie) que le
corps engendr\'e par $P_{n,c,4}(x)$ l'est aussi par
${\tilde{P}}_{n,c,4}(x) = x^4-2x^3+(1-n)x^2+nx-c$. Le calcul montre,
pour une sp\'ecialisation de $(n,c)= (\frac{t^2+32}{2t^2},
\frac{3t^4-1024}{16t^4})$, que  
$${\mbox{\rm{R\'esultant}}}({\tilde{P}}_{n,c,4}(x), X-(\frac{t}{2}x^2 -
  \frac{t^2+32}{8t}),x) = X^4 - tX^3 - 6X^2 + tX + 1,$$
et l'on retrouve ainsi la famille quartique cyclique de Gras.

\end{document}